\documentclass[leqno,12pt]{amsart} 
\usepackage{amsmath, amssymb}
\usepackage{amsthm} 
\usepackage[all]{xy}
\usepackage{pstricks, pst-tree}
\usepackage{enumerate}
\theoremstyle{plain} 
\newtheorem{theorem}{\indent\sc Theorem}[section]
\newtheorem{lemma}[theorem]{\indent\sc Lemma}
\newtheorem{corollary}[theorem]{\indent\sc Corollary}
\newtheorem{proposition}[theorem]{\indent\sc Proposition}

\theoremstyle{definition} 
\newtheorem{definition}[theorem]{\indent\sc Definition}
\newtheorem{remark}[theorem]{\indent\sc Remark}
\newtheorem{example}[theorem]{\indent\sc Example}

\newcommand{\N}{\mathbb{N}}
\newcommand{\Z}{\mathbb{Z}}
\newcommand{\R}{\mathbb{R}}
\newcommand{\C}{\mathbb{C}}

\newcommand{\abs}[1]{\left\lvert#1\right\rvert}
\newcommand{\norm}[1]{\lVert#1\rVert}
\newcommand{\Chig}{C_{h_L}}
\newcommand{\CHig}{\overline{\Chig}}

\newcommand{\cx}{\times_{\mathrm{cone}}}

\makeatletter
   
   \@addtoreset{equation}{section}
\makeatother
\title{\uppercase{Sublinear Higson corornae of Euclidean cones}}
\author{
\textsc{Tomohiro Fukaya} 
}
\date{} 
\begin{document}

\maketitle

\footnote{ 
2000 \textit{Mathematics Subject Classification}.
Primary 51F99; Secondary 46L45.
}
\footnote{ 
\textit{Key words and phrases}. 
Coarse geometry, Higson corona, homotopy, commutative $C^*$-algebra
}
\footnote{ 

}

\begin{abstract}
 The aim of this paper is to introduce the sublinear Higson corona 
 and show that the sublinear Higson corona of Euclidean cone of $P$ and $X$ 
 is decomposed into the product of $P$ and that of $X$. Here $P$ is a
 compact metric space and $X$ is unbounded proper metric space.
 For example, the sublinear Higson corona of  $n$-dimensional Euclidean
 space is homeomorphic to  the product of $(n-1)$-dimensional sphere and 
 that of natural numbers.
\end{abstract}

\section{introduction}
In the coarse geometry, it is very important to study boundaries 
of metric spaces. If the metric space $X$ is hyperbolic in
the sense of Gromov, we can define the ideal boundary and
construct a natural metric on it. There are a lot of studies
of relations between the coarse geometry of $X$ and topology of the
ideal boundary of $X$. 

For general proper metric spaces, we can define the boundary,
which is called Higson corona. It is a contravariant functor from the
category of coarse spaces into that of compact Hausdorff spaces.
Since Higson coronae are never metrizable, their topology seem extremely
complicated. However, there are several studies of the topology of
Higson coronae. For example, if
the asymptotic dimension of the metric space is finite, then it is equal
to the covering dimension of its Higson corona \cite{MR1840358}. 
Dranishnikov and Ferry showed that $n$-dimensional Euclidean space and
Hyperbolic space of the same dimension have different $n$-th \v{C}ech
cohomologies~\cite{MR1490028}. 
The author studied the
group action and fixed point theorem on Higson
coronae~\cite{fukaya-coarse-dynamics}.

In this paper, we introduce the sublinear Higson corona $\nu_L X$ for
coarse space $X$ and show that sublinear Higson corona of Euclidean cone
$P\cx X$ is decomposed into the product of $P$ and $\nu_L X$
(Theorem~\ref{main_theorem}). Here $P$ is a  compact metric space
and $P\cx X$ is defined in Section~\ref{sec:euclidean-cone}.
For example, $\nu_L \R^n$ is
homeomorphic to $S^{n-1}\times \nu_L \N$. As an application, we 
consider the linear map $T\colon \R^n\rightarrow \R^n$ of a positive
determinant. We show that the induced map on $\nu_L \R^n$
is homotopic to the identity. We remark that this statement does not
hold on the Higson corona $\nu \R^n$.

The organization of this paper is as follows: 
In section~\ref{sec:line-coarse-struct}, we introduce the 
coarse category. Our definition is slightly different from usual one.
In section~\ref{sec:higs-comp-line}, we define the
sublinear Higson compactification and study its functorial
properties. In section~\ref{sec:euclidean-cone}, we define the Euclidean cone
and study the decomposition of  the sublinear Higson corona of it.
In section~\ref{sec:applications}, we give an application. 
In Appendix, we give a proof of Lemma~\ref{lem:coarse_and_continuous}.

\subsubsection*{Acknowledgment}
The author would like to thank Makoto Yamashita for many helpful
discussions. The author was supported by 
Grant-in-Aid for JSPS Fellow (19$\cdot$3177) from Japan Society for the
Promotion of Science.

\section{coarse category}
\label{sec:line-coarse-struct}
A metric space $X$ is {\itshape proper} if every closed bounded set in $X$
is compact. For a positive number $C$, we say that $X$ is 
{\itshape $C$-quasi-geodesic space} if for any $x,x'\in X$, there exists
a map $f \colon [0, d(x,y)] \rightarrow X$ such that 
$(1/C) \abs{a-b} -C \leq d(f(a),f(b)) \leq C\abs{a-b} +C$ for any 
$a,b \in [0, d(x,y)]$. We choose a base point $e\in
X$ and define $\abs{x} := d(e,x)$ for $x\in X$. 
We say that $X$ is 
{\itshape coarse space} if $X$ is proper and $C$-quasi-geodesic space for
some $C$ with a base point $e$. 
\begin{definition}
 Let $X$ and $Y$ be coarse spaces and let $f\colon X\rightarrow Y$
 be a map (not necessarily continuous). We say that the map $f$ is 
 {\itshape coarse} if there exists a positive constant $A$ such
 that the following two conditions are satisfied.
 \begin{enumerate}
 \item $\abs{f(x)} \geq \abs{x}/A-A$ for every $x\in X$.
 \item $d(f(x),f(x')) \leq Ad(x,x') +A$ for every $x,x' \in X$. 
 \end{enumerate}
\end{definition}
Let $f,g \colon X\rightarrow Y$ be maps. We define that $f$ is 
{\itshape sublinearly close} to $g$ if for any $\epsilon > 0$,
there exists a positive constant $C_\epsilon$ such that
$d(f(x),g(x)) \leq \epsilon \abs{x} + C_\epsilon$ for all $x\in X$.

We define that metric spaces
$X$ and $Y$ are {\itshape coarsely equivalent} if there exist
coarse maps $f\colon X\rightarrow Y$ and $g\colon Y\rightarrow
X$ such that $g\circ f$ and $f\circ g$ are sublinearly close to the identity
maps of $X$ and $Y$, respectively. The category of coarse spaces
consists of coarse spaces and coarse maps. 

\begin{remark}
 Our definition of the coarse equivalence is different from that of for metric
 spaces with usual bounded coarse structure. Compare with Definition 1.8
 and Example 2.5 of \cite{MR2007488}. Our definition is related to the
 sublinear coarse structure defined in \cite{MR2294641}.
\end{remark}

\section{Sublinear Higson corona}
\label{sec:higs-comp-line}
All arguments in this section is based on that in Section 2.3 of
\cite{MR2007488}.
Let $X$ be a coarse space and $\varphi\colon X \rightarrow \C$ 
be a bounded continuous map. 
We say that a $\varphi$ is sublinear Higson function on $X$
if there exists a constant $C_\varphi$ such that for all $R> 0$ and
$x,x' \in X\setminus B(R)$, we have
\begin{eqnarray}
\label{eq:Chig}
  \abs{\varphi(x) - \varphi(x')} < 
 \frac{C_\varphi d(x,x')}R.
\end{eqnarray}
Here $B(R)$ denotes the ball of radius $R$ centered at the base point $e$.
We let $\Chig(X)$ denote the space of all bounded continuous sublinear
Higson functions on $X$.
Then $\Chig(X)$ is a sub $*$-algebra of $C_b(X)$ and its closure,
denoted by $\CHig(X)$, is a unital $C^*$-algebra.

By the Gelfand-Naimark theorem, $\CHig(X)$ is the algebra of continuous
functions on some compactification of $X$.
\begin{definition}
 The compactification $h_L X$ of $X$ characterized by the property
 $C(h_L X) = \CHig(X)$ is called the 
{\itshape sublinear Higson compactification}. Its boundary 
$h_L X \setminus X$ is denoted $\nu_L X$, and is called the 
{\itshape sublinear Higson corona} of $X$.
\end{definition}
We let $C_0(X)$ denote the algebra of continuous functions on $X$ which
banish at infinity. Then we have $C(\nu_L X) \cong \CHig(X)/C_0(X)$.
\begin{example}
 We define an embedding 
$i \colon \R^2=\C \hookrightarrow D^2= \{z\in \C\colon \abs{z}\leq 1\}$
 by sending $z$ to $z/(1+\abs{z})$. For any
 $\varphi \in C(D^2)$, the composite $\varphi\circ i$ belongs to 
$\CHig(\R^2)$. Thus $i$ induces the injection 
$i^*\colon C(D^2) \hookrightarrow \CHig(R^2)$, it follows that
$i$ extends to the surjection
\begin{eqnarray*}
 h_L i \colon h_L \R^2 \rightarrow D^2.
\end{eqnarray*} 
\end{example}
\begin{proposition}
 As is Higson corona, the sublinear Higson corona of an unbounded coarse
 space is never second countable and its cardinal number is
 greater than or equal to $2^{2^{\aleph_0}}$.
\end{proposition}
\begin{proof}
 For second countability, it is enough to show that $\CHig(X)$ is not
 separable. 
 We can
 choose a sequence $\{x_n\}$ such that $\abs{x_n} >
 2\abs{x_{n-1}}$ for all $n\geq 0$. 
 We define a continuous map $\varphi_n \colon X \rightarrow \C$ as follows:
\[
 \varphi_n(x) =
 \begin{cases}
  1-\frac{4d(x,x_n)}{\abs{x_n}}& 
    \text{ if } d(x,x_n) \leq \frac{\abs{x_n}}{4}, \\
  0 & \text{ otherwise.}
 \end{cases}
\]
 For any map $P \colon \N \rightarrow \{0,1\}$, we define the continuous map
 $\psi\colon X\rightarrow \C$ by 
 $\psi_P(x) = \sum_{n\in \N} P(n)\varphi_n(x)$. Thus we obtain an
 uncountably many family of sublinear Higson functions such that
 $\norm{\psi_P - \psi_{P'}} = 1$ for 
 any pair $(\psi_P, \psi_{P'})$ of distinct $P,P'\in \{0,1\}^\N$. This
 shows that $\CHig(X)$ is not separable.
 For any $\psi \in C_b(\{x_n\}) = l^\infty (\{x_n\})$,
 an extension $\tilde{\psi}\in \CHig(X)$ is given by
 $\tilde{\psi}(x) = \sum_{n\in \N}\psi(x_n)\varphi_n(x)$.
 This means the restriction 
 $\CHig(X) \rightarrow C_b(\{x_n\})$ is surjective and thus the inclusion
 $\{x_n\} \hookrightarrow X$ extends to the embedding 
 $\beta\{x_n\} \rightarrow h_L X$. 
 Since $\beta\{x_n\}$ is homeomorphic to $\beta\N$, the cardinal number
 of $h_L X$ is greater than or equal to that of $\beta\N$,
 that is, $2^{2^{\aleph_0}}$.
\end{proof}
Let $f\colon X\rightarrow Y$ be a continuous coarse map. Then $f$
extends to the continuous map $h_Lf\colon h_LX\rightarrow h_LY$. 
The restriction to $\nu_L X$ is denoted $\nu_L f:= h_Lf|_{\nu_L X}$. 
However, even if $f$ is coarse
but not continuous, we can define the map 
$\nu_L f\colon \nu_L X \rightarrow \nu_L Y$ as follows.

Let $X$ be a coarse space.
Let $B_{h_L}(X)$ denote the algebra of all bounded (not necessary
continuous) functions 
$\varphi\colon X \rightarrow \C$ which satisfy the following condition:
 There exists a constant $C_\varphi$ such that for all $R> 0$ and
 $x,x' \in X\setminus B(R)$,
 we have 
\[
 \abs{\varphi(x) - \varphi(x')} < 
 \frac{C_\varphi d(x,x') +C_\varphi}{R}.
\]
Let $\overline{B_{h_L}}(X)$ be a closure of $B_{h_L}(X)$ in the
algebra $B(X)$ of all bounded functions $X \rightarrow \C$ and 
$B_0(X)$ be the ideal of all bounded functions that tend to $0$
at infinity.
\begin{definition}
 A coarse space $X$ is said to be of {\itshape bounded geometry} if
 there exists a
 uniformly bounded cover $\mathcal{U}=\{U_\alpha\}$ of $X$
 with positive Lebesgue number and of finite degree. That is, there
 exist constants $L, d$ 
 and $N$ such that Lebesgue number of $\mathcal{U}$ is $L$, 
 the diameter of all $U_\alpha \in U$ is bounded by
 $d$ and no more than $N$ member of $\mathcal{U}$ has non-empty
 intersection.
\end{definition}
\begin{lemma}
\label{lem:coarse_and_continuous}
 Let $X$ be a coarse space of bounded geometry.
 Then
\begin{enumerate}[$(a)$]
 \item $C_0(X) = \CHig(X) \cap B_0(X)$.
\label{item:a}
 \item $\overline{B_{h_L}}(X) = \CHig(X) + B_0(X)$.
\label{item:b}
\end{enumerate}
\end{lemma}
The part~(\ref{item:a}) is obvious. For the proof of the part~
(\ref{item:b}), see Appendix. 

It follows from the second Isomorphism theorem that
\[
 C(\nu_L X) = \frac{\CHig(X)}{C_0(X)}=
 \frac{\CHig(X)}{\CHig(X)\cap B_{0}(X)}=
 \frac{B_0(X) + \CHig(X)}{B_{0}(X)}=
 \frac{\overline{B_{h_L}}(X)}{B_0(X)}.
\]

\begin{proposition}
\label{prop:higs-comp-line-1}
 Let $X$ and $Y$ be coarse spaces of bounded geometry.
 A coarse map $f \colon X \rightarrow Y$ extends to a
 continuous map $\nu_L f \colon \nu_L X \rightarrow \nu_L Y$.
If $f,g\colon X\rightarrow Y$ are sublinearly close, then $\nu_L f=\nu_L g$.
\end{proposition}
\begin{proof}
 A coarse map $f\colon X\rightarrow Y$ induces
 $*$-homomorphism 
 $f^*\colon \overline{B_{h_L}}(Y)\rightarrow \overline{B_{h_L}}(X)$ and 
 $f^*\colon B_{0}(Y)\rightarrow B_{0}(X)$. If $f$
 is sublinearly close to $g$ then $f^*-g^*$ maps
 $\overline{B_{h_L}}(Y)$ to $B_0(X)$.
\end{proof}
\begin{corollary}
\label{cor:higs-comp-line-2}
 Coarsely equivalent spaces of bounded geometry
 have homeomorphic sublinear Higson corona.
\end{corollary}

\begin{proposition}
 Sublinear Higson corona is a faithful functor from the category of
 coarse spaces of bounded geometry to that of compact Hausdorff spaces. 
 That is, for any pair of coarse maps $f,g$, the equality 
$\nu_L f = \nu_L g$ implies that $f$ is sublinearly close to $g$.
\end{proposition}
\begin{proof}
Let $f,g \colon X\rightarrow Y$ be coarse maps. 
We suppose that $f$ is not sublinearly close to $g$. There
 exist a  constant $C$ and a sequence $\{x_n\}$ such that 
$d(f(x_n),g(x_n))\geq C\abs{x}$. We can construct a sublinear Higson
 function $\varphi$ on $Y$ such that $\varphi(f(x_n)) = 1$ and 
 $\varphi(g(x_n)) = 0$ for all $n \geq 0$, therefore $\nu_L f \neq \nu_L g$.
\end{proof}

\begin{remark}
 Dranishnikov and Smith introduced the {\itshape sublinear coarse
 structure} \cite{MR2294641} and showed that sublinear Higson corona of proper
 metric space is 
 homeomorphic to its Higson corona with respect to the sublinear coarse
 structure (Theorem 2.11 of \cite{MR2294641}).
\end{remark}

\section{Euclidean cone}
\label{sec:euclidean-cone}
Let $X$ be a coarse space and $P$ be a compact path metric space. 
We define
the Euclidean (punctured) cone metric on $P\times X$ as follows. A {\itshape path}
$\gamma$ is a continuous map from an interval 
$I\subset [0,\infty)$ to $P\times X$. 
We define the {\itshape length} of $\gamma$ to be
\[
 \sup\left\{\sum_{j=o}^{n-1}d(x_j,x_{j+1}) +
 \max\{1, \abs{x_j},\abs{x_{j+1}}\}d(p_j,p_{j+1})\right\},
\]
where the supremum is taken over all finite sequences
$(p_j,x_j)_{j=0}^n$ of points on the path $\gamma$ with 
$(p_0,x_0)$ and $(p_n,x_n)$ being the two endpoints. We make $P\times X$
into a path metric space by defining the distance $d$ between two points
to be the infimum of the length of all possible paths joining them.
The {\itshape Euclidean cone} of $P$ and $X$, denoted by $P\cx X$, is
the metric space $P\times X$ equipped with this metric. 
\begin{example}
 The Euclidean cone $S^{n-1} \cx \N$ is coarsely equivalent to $\R^n$.
\end{example}

We define the map 
$\Lambda\colon \Chig(P\cx X) \rightarrow C(P\times \CHig(X))$ by sending
$\varphi \in \Chig(Y)$ to the map which sends $p\in P$ to 
$\varphi_p$:
\begin{eqnarray*}
 \Lambda \colon 
 \Chig(P\cx X))\ni \varphi \mapsto (p \mapsto \varphi_p) \in C(P,\CHig(X)).
\end{eqnarray*}
Here $\varphi_p$ denotes the map $x\mapsto \varphi(p,x)$ and
$C(P,\CHig(X))$ denotes the $C^*$-algebra of continuous maps from
$P$ to $\CHig(X)$.
\begin{proposition}
\label{prop:lambda}
 $\Lambda$ is well-defined and extends to 
 a $*$-monomorphism of $C^*$-algebras:
 \[
  \Lambda\colon \CHig(P\cx X)\rightarrow C(P,\CHig(X)).
 \] 
\end{proposition}
\begin{proof}
 Let $\varphi \in \Chig(P\cx X)$.
 It is clear that $\varphi_p$ belongs to $\CHig(X)$ for all 
 $p\in P$.
  We show that the map
 $\Lambda(\varphi)\colon P\rightarrow \CHig(X)$ is continuous.
 Let $p,p'\in P$. For all $x\in X$, 
\begin{eqnarray*}
 \abs{\Lambda(\varphi)(p)(x)-\Lambda(\varphi)(p')(x)} &=& 
  \abs{\varphi(p,x) - \varphi(p',x)}\\
 &\leq& \frac{C_\varphi}{\abs{x}} d((p,x),(p',x))\\
 &\leq& C_\varphi d(p,p') .
\end{eqnarray*}
 Thus we have $\norm{\Lambda(\varphi)(p)-\Lambda(\varphi)(p')}
 \leq C_\varphi d(p,p')$ and therefore 
 $\Lambda(\varphi) \in C(P,\CHig(X))$.
 Since $\Lambda$ is a Lipschitz map with Lipschitz constant less than
 or equal to 1, we see that $\Lambda$ extends to $\CHig(Y)$. 
 It is clear that $\Lambda$ is a $*$-monomorphism.
\end{proof}
\begin{proposition}
\label{prop:Omega}
$\Omega\colon C(P)\otimes \CHig(X) \rightarrow \CHig(P\cx X)
 : \varphi \otimes \psi \mapsto \varphi\cdot \psi$ is well-defined. 
\end{proposition}
\begin{proof}
 Let $\varphi \otimes \psi \in C(P)\otimes \Chig(X)$. 
 We assume that $\varphi$ is a Lipschitz map with a Lipschitz constant 
 $C_\varphi$. Let $R> 0$. For any 
 $(p,x),(p',x') \in P\cx X$ with $\abs{x},\abs{x'}>R$, we have 
 \[
  d(x,x') + \abs{x}d(p,p') \leq d((p,x),(p',x')).
 \] 
 Thus 
\begin{eqnarray*}
  \abs{\varphi(p) \psi(x) - \varphi(p') \psi(x')} &\leq& 
 \abs{\varphi(p)\psi(x) - \varphi(p') \psi(x)} 
 +\abs{\varphi(p') \psi(x) - \varphi(p')\psi(x')}\\
 &\leq& C_\varphi\norm{\psi}d(p,p')
  + \frac{C_\psi\norm{\varphi}d(x,x')}{R}\\
 &\leq& \frac{C_\varphi\norm{\psi} + C_\psi\norm{\varphi}}{R}d((p,x),(p',x')).
\end{eqnarray*}
 It follows that $\Omega(\varphi \otimes \psi)$ belongs to
 $\CHig(P\cx X)$. Since the set of Lipschitz maps is dense in $C(P)$,
 we have the desired consequence.
\end{proof}

To show that $\Omega$ gives the inverse of $\Lambda$, we need the
following well-known fact.
\begin{lemma}
\label{lem:well-known-fact}
 Let $P$ be a compact metric space and $A$ be a commutative
 $C^*$-algebra. Then $C(P)\otimes A \cong C(P,A)$.
\end{lemma}
\begin{proof}
 We can construct a family $\{\mathcal{U}^n\}_{n\in \N}$ of finite
 covers of $P$ such that, the diameter of each member
 $U_i^n$ of $\mathcal{U}^n$ is less than $1/n$. We choose points $p_i^n\in
 U_i^n$ for each $n\in N$. Let $\{h_i^n\}_i$ be a
 partition of unity subordinate to $\mathcal{U}^n$. 
 We define $\Psi\colon C(P)\otimes A \rightarrow  C(P,A)$
 by $\Psi(\varphi\otimes a)(p) = \varphi(p)a$
 for $\varphi \in C(P),\, a\in A$ and $p\in P$. Clearly $\Psi$ is 
 injective. We show that $\Psi$ is surjective. Let $\psi \in C(P,A)$.
 Set $\psi_n = \sum_{i}h_i^n\otimes \psi(p_i^n)\in C(P)\otimes A$.
 Since $\psi$ is uniformly continuous, 
 $\norm{\psi - \Psi(\psi_n)}$ tends to $0$ as $n$ goes to infinity.
 Thus $\Psi$ is surjective.
\end{proof}
\begin{theorem}
\label{main_theorem}
 The sublinear Higson compactification of the Euclidean cone $P\cx X$ is
 homeomorphic to the product $P\times h_LX$. Especially 
$\nu_L(P\cx X) = P\times \nu_LX$.
\end{theorem}
\begin{proof}
 By Proposition~\ref{prop:lambda},~\ref{prop:Omega} and
 Lemma~\ref{lem:well-known-fact}, 
 $\CHig(P\cx X)\cong C(P)\otimes \CHig(X)$.
\end{proof}
\begin{example}
 $\nu_L(\R^n) = S^{n-1} \times \nu_L\N$.
\end{example}
\section{Applications}
\label{sec:applications}
\begin{definition}
 Let $f,g\colon X\rightarrow Y$ be coarse maps. We say that $f$
 is {\itshape cone-homotopic} to $g$ if there exists a coarse map
$H\colon [0,1]\cx X\rightarrow Y$  such that $f=H_0$ and $g=H_1$. Such
 an $H$ is called {\itshape cone homotopy} between $f$ and $g$.
\end{definition}
\begin{theorem}
 If $f$ is cone-homotopic to $g$, then the induced map 
 $\nu_L f$ is homotopic to $\nu_L g$.
\end{theorem}
\begin{proof}
 $H$ induces a continuous map 
 $\nu_LH\colon [0,1]\times \nu_LX\rightarrow \nu_LY$
 such that $\nu_LH(0,x) = \nu_Lf(x)$ and $\nu_LH(1,x) = \nu_Lg(x)$ for
 all $x\in X$.
\end{proof}
\begin{example}
 Let $T$ be an $n$ by $n$ integer matrix with a positive determinant.
 Then linear map $T\colon \Z^n \rightarrow \Z^n$ 
 is a coarse map. The induced map 
 $\nu_L T\colon \nu_L \Z^n \rightarrow \nu_L \Z^n$ is
 homotopic to the identity map $\mathrm{id}_{\nu_L \Z^n}$.
\end{example}
We remark that since $T$ is not close to the 
identity $I_n$, the induced map $\nu_L T$ is different from the
identity  
$\mathrm{id}_{\nu_L \Z^n}$.
\begin{proof}
 Since $\nu_L \Z^n$ is homeomorphic to $\nu_L \R^n$, 
 it is enough to show that the map 
$T\colon \R^n \rightarrow \R^n$ is cone-homotopic to the identity map. 
Since $T$ has a positive determinant, we can choose a continuous path 
$\Theta\colon [0,1] \rightarrow GL_+(n,\R)=\{A\in GL(n,\R): \det A >0\}$
such that $\Theta(0)  = T$ and $\Theta(1) = I_n$. A map 
$H(x,t) = \Theta(t)x$ is a cone homotopy  between $T$ and the identity $I_n$.
\end{proof}
The same statement for the Higson corona of $\Z^n$ does not hold. This
is pointed out by Makoto Yamashita. Here we recall the definition of the
Higson corona. Let $X$ be a coarse space.
We define a $C^*$-algebra by
\[
 C_h(X) = \{f\in C_b(X) : 
\lim_{\abs{x}\rightarrow \infty} \text{diam}(f(B(x,r))) = 0 \quad
{}^\forall r>0\}.
\]
Here $B(x,r)$ denotes the $r$-ball centered at $x\in X$.
The Higson compactification of $X$, denoted $hX$, is the
compactification characterized by $C_h(X) = C(hX)$. The Higson corona,
denoted $\nu X $, is defined by $hX \setminus X$. A continuous coarse
map $f\colon X\rightarrow Y$ extends to $hf\colon hX \rightarrow hY$. 
The restriction to $\nu X$ is denoted $\nu f := hf|_{\nu X}$.
The following argument is based on that in Section 3 of \cite{MR1369755}. 
\begin{lemma}
\label{lem:homotopy-l-higson}
 For $s>0$, we define a map  $\phi_s \colon \R_+ \rightarrow \R$
 by $x\mapsto \sqrt{sx}$. 
 Let $e\colon \R\rightarrow S^1$ be a covering $e(x) = \exp(\sqrt{-1}x)$.
 The composite $e\circ \phi_s$ extends $h(e\circ \phi_s) \colon
 h\R_+\rightarrow S^1$. If $h(e\circ \phi_s)$ is homotopic to 
 $h(e\circ \phi_t)$ for some $t>0$, then $s=t$.
\end{lemma}
\begin{proof}
 Since $\text{diam}(\phi_s(B(x,r)))\rightarrow 0$ as 
 $\abs{x}\rightarrow \infty$ for any $r>0$, we have a $*$-homomorphism 
 $(e\circ \phi_s)^*\colon C(S^1)\rightarrow C_h(\R_+)$. Thus we have the
 extension $h(e\circ \phi_s)\colon h\R_+\rightarrow S^1$.
 We suppose that $h(e\circ \phi_s)$ is homotopic to 
 $h(e\circ \phi_t)$. We define $\psi_{s,t}\colon h\R_+\rightarrow S^1$ by 
\[
 \psi_{s,t}(x) = \frac{h(e\circ \phi_s)(x)}{h(e\circ \phi_t)(x)}.
\]
 Since $\psi_{s,t}$ is null-homotopic, there
 exists a lift $\psi_{s,t}^{\sim}\colon hX\rightarrow \R$ such that
 $e\circ \psi_{s,t}^\sim = \psi_{s,t}$. The restriction
 $\psi_{s,t}^\sim|_X$ and $\phi_s - \phi_t$ are both lifts of
 $\psi_{s,t}|_X$. The image of $\psi_{s,t}^\sim$
 must be bounded since $hX$ is compact. It follows that the image of 
 $\phi_s - \phi_t$ is also bounded. Thus we have $s=t$.
\end{proof}
\begin{proposition}
For $s>0$, we define a map $f_s\colon \R_+\rightarrow \R_+$ by 
$f_s(x) =sx$. If the induced map 
$\nu f_s \colon \nu \R_+ \rightarrow \nu \R_+$ is homotopic to $\nu f_t$ for
some $t>0$. Then $s =t$.
\end{proposition}

\begin{proof}
We suppose the induced map 
$\nu f_s \colon \nu \R_+ \rightarrow \nu \R_+$ is homotopic to $\nu f_t$ for
some $t>0$. Let $H\colon  [0,1]\times \nu \R_+ \rightarrow \nu S^1$ be
 a homotopy between $h(e \circ \phi_1)\circ \nu f_s$ and 
$h(e \circ \phi_1)\circ \nu f_t$.
Set $A = [0,1] \times \nu \R^n \cup  \{0,1\} \times h\R^n
\subset  \{0,1\} \times \R^n$.
We define a map $H'\colon A\rightarrow S^1$ by 
\begin{eqnarray*}
 H'(u,x)=
 \begin{cases}
 H(u,x) & \text{if } 0<u<1,\\
 h(e\circ \phi_s)(x) & \text{if } u=0,\\
 h(e\circ \phi_t)(x) & \text{if } u=1.
 \end{cases}
\end{eqnarray*} 
 Here we remark $ h(e\circ \phi_s) = h(e\circ \phi_1)\circ hf_s$.
 Since $S^1$ is ANR, there exist a neighborhood $U$ of $\nu X$ in $hX$
 and a extension 
 $H'' \colon U\times[0,1] \rightarrow S^1$ such that 
 $H''|_{A\cap U} = H'|_{A\cap U}$.
 Then we have $h(e\circ \phi_s)|_U$ 
 is homotopic to $h(e\circ \phi_t)|_U$.
 By Lemma~\ref{lem:homotopy-l-higson}, we have $s=t$. 
\end{proof}

\begin{corollary}
 Let $T\colon \Z^n \rightarrow \Z^n$ be a linear map with an eigenvalue
 $\lambda \neq 1$. Then the induced map 
$\nu T\colon \nu \Z^n\rightarrow \nu \Z^n$ is not homotopic to the
 identity $\mathrm{id}_{\nu \Z^n}$.
\end{corollary}
\section{Appendix--proof of Lemma~\ref{lem:coarse_and_continuous}}
\label{sec:appendix}
The proof is based on the proof of Lemma 2.40 of \cite{MR2007488}.
By the assumption of $X$, there exists a cover
$\mathcal{U} = \{U_\alpha\}$ of $X$ and 
constants $L,d,N$ such that Lebesgue number of $\mathcal{U}$ is $L$,
the diameter of any $U_\alpha$ is less than $d$ and no more than $N$
members of $\mathcal{U}$ has non-empty intersection.
Then we can construct a partition of unity $\pi_\alpha$ subordinate to 
$\mathcal{U}$ all of whose constituent functions are 
$D$-Lipschitz for some constant $D=D(L,d,N)$.
(See the proof of Theorem 9.9 of \cite{MR2007488}.)

We choose a point $x_\alpha \in U_\alpha$ for each $\alpha$.
Now we let $f \in B_{h_L}(X)$ and $C_f$ be a constant which
appears in the definition of $B_{h_L}(X)$. We define 
\[
 g(x) := \sum_{\alpha}\pi_\alpha(x)f(x_\alpha).
\]
The function $g$ is continuous and bounded. For all $x\in X$, we have
\[
 f(x) - g(x) = \sum_{\alpha}\pi_\alpha(x)(f(x) - f(x_\alpha))
\]
and $d(x,x_\alpha) < d$ whenever $\pi_\alpha(x) \neq 0$. Thus we have 
$f-g \in B_0(X)$. Next we will show that $g$ satisfies (\ref{eq:Chig}).
We assume that $X$ is a $C_X$-quasi-geodesic for some positive constant
$C_X$. 
Let $R > 2d$ and $x,x' \in X\setminus B(R)$ such that 
$d(x,x') \leq C_X$.
Set $I_{+} := \{\alpha \colon \pi_\alpha(x) - \pi_\alpha(x') > 0\}$
and $I_{-} := \{\alpha \colon \pi_\alpha(x) - \pi_\alpha(x') < 0\}$. 
Set $t:= \sum_{\alpha \in I_{+}}(\pi_\alpha(x) - \pi_\alpha(x')) = 
- \sum_{\alpha \in I_{-}}(\pi_\alpha(x) - \pi_\alpha(x'))$. Since each 
$\pi_\alpha$ is $D$-Lipschitz, we have $t\leq 2NDd(x,x')$.
Set 
$f_{\max} := \max\{f(x_\alpha)\colon \alpha \in I_{+}\}$ and 
$f_{\min} := \min\{f(x_\alpha)\colon \alpha \in I_{-}\}$. 
For any $\alpha, \alpha' \in I_{+}\cup I_{-}$, 
we have $x_{\alpha},x_{\alpha'} \in X \setminus B(R-d)$ and 
$d(x_\alpha, x_{\alpha'}) < C_X + 2d$. It follows that 
$f_{\max} - f_{\min} \leq C_f(C_X + 2d)/(R-d) 
< 2C_f(C_X + 2d)/R$.
We can assume that $g(x) \geq g(x')$. Then we have
\begin{eqnarray*}
 \abs{g(x) - g(x')} &=& g(x) - g(x') \\
 &=& \sum_{\alpha \in I_{+}}(\pi_\alpha(x) - \pi_\alpha(x'))f(x_\alpha)
 + \sum_{\alpha \in I_{-}}(\pi_\alpha(x) - \pi_\alpha(x'))f(x_\alpha)\\
 &\leq& \sum_{\alpha \in I_{+}}(\pi_\alpha(x) - \pi_\alpha(x'))f_{\max}
   + \sum_{\alpha \in I_{-}}(\pi_\alpha(x) - \pi_\alpha(x'))f_{\min}\\
 &=& t(f_{\max} - f_{\min}) \\
 &\leq& \frac{4NDC_f(C_X+2d)}{R}d(x,x').
\end{eqnarray*}
This shows that $g\in \CHig(X)$ and completes the proof of
Lemma~\ref{lem:coarse_and_continuous}.

\bibliographystyle{amsplain}
\bibliography{/DirUsers/tomo_xi/Library/tex/books,/DirUsers/tomo_xi/Library/tex/math}

\bigskip
\address{ Tomohiro Fukaya \endgraf
Department of Mathematics, Kyoto University, Kyoto 606-8502, Japan}

\textit{E-mail address}: \texttt{tomo\_xi@math.kyoto-u.ac.jp}
\end{document}